\documentclass[12pt, twoside]{amsart}
\usepackage[english]{babel}
\usepackage{amssymb}
\usepackage{amsmath}
\usepackage{latexsym}
\usepackage{amsthm}
\usepackage{mathrsfs}
\usepackage{enumerate}
\theoremstyle{plain}
\numberwithin{equation}{section}
\addtocounter{page}{1}
\newtheorem{theorem}{Theorem}[section]
\newtheorem{remark}{Remark}[section]
\newtheorem{proposition}{Proposition}[section]

\newtheorem{definition}{Definition}[section]
\newtheorem{example}{Example}[section]

\newcommand{\n}{|\!|}
\title{Local solvability of a class of degenerate second order operators
\\
}
\author{Serena Federico
}
\address{DIPARTIMENTO DI MATEMATICA, UNIVERSIT\`{A} DI BOLOGNA, PIAZZA DI PORTA S. DONATO 5, 40126, BOLOGNA, ITALIA.
}
\email{serena.federico2@unibo.it
}

\begin{document}

\begin{abstract}
{In this paper we will first present some results about the local solvability property of a class of degenerate second order linear partial differential operators with smooth coefficients.
 The class under consideration (which in turn is a generalization of the Kannai operator) exhibits a degeneracy due to the interplay between the singularity associated with the characteristic set of a system of vector fields and the vanishing of a function. Afterward we will also discuss some local solvability results for two classes of degenerate second order linear partial differential operators with non-smooth coefficients which are a variation of the main class presented above.
}


\end{abstract}
\maketitle
\section{
Introduction}
This article deals with the local solvability problem of some classes of second order linear degenerate partial differential operators with smooth and non-smooth coefficients. The main class, which is the one with smooth coefficients, is given by
\begin{equation}
\label{Pintro}
P=\sum_{j=1}^NX_j^*f^pX_j+iX_0+a_0,
\end{equation}
where $p\geq 1$ is an integer, $X_j(x,D)$, $0\leq j\leq N$ ($D=-i\partial$), are \textit{homogeneous first order differential operators} 
(i.e. with no lower order terms) with smooth coefficients on an open set $\Omega\subset\mathbb{R}^n$
and with a \textit{real} principal symbol (in other words, the $iX_j(x,D)$ are \textit{real} vector fields), $a_0$ is a smooth possibly complex-valued function, and $f\colon\Omega\longrightarrow\mathbb{R}$ a smooth function with $f^{-1}(0)\not=\emptyset$ and $df\bigl|_{f^{-1}(0)}\not=0$. We also suppose that $iX_0f>0$ on $S=f^{-1}(0)$, that is, the real vector field $iX_0$ is transverse to $S$ and can only be zero at some point of $\Omega\setminus S$. 

Due to the presence of the function $f^p$ (which is vanishing on $S$)  in the leading part of the operator, the class \eqref{Pintro} exhibits a degeneracy on $S$ due to the interplay between the characteristic set of the system of vector fields $\{iX_j\}_{j=1,...,N}$ and the vanishing of the function $f$, therefore \eqref{Pintro} contains degenerate multiple characteristics partial differential operators whose degree of degeneracy depends on the parity of the exponent $p$.  
Moreover, since non-vanishing conditions are not imposed on the vector fields  $iX_j$, $1\leq j\leq N$, appearing in the second order part of \eqref{Pintro}, we have that the operators in \eqref{Pintro} can be degenerate far from $S$ as well. 

By virtue of the properties mentioned above the local solvability of the class under consideration is not guaranteed, and, of course, the local solvability of $P$ at points of $S$ is more difficult to analyze. On $S$ the most interesting situation occurs when $p$ is an odd integer since we have that the principal symbol of the operator changes sign in the neighborhood of each point of the set $\pi^{-1}(S)$, where $\pi^{-1}(S)$ denotes the cotangent fiber of $S$.

The changing of sign of the principal symbol is a property that can obstacle the local solvability of an operator at points where the property holds. This is indeed what happens for the celebrated Kannai operator  $P=x_1\sum_{j=2}^nD_{j}^2-iD_1$ (having multiple characteristics), which is hypoelliptic over the whole $\mathbb{R}^n$ but is not locally solvable at the set $S=\{x\in\mathbb{R}^n;x_1=0\}$ around which its principal symbol changes sign (see \cite{Ka}).

This example also highlights the relation between hypoellipticity and local solvability. It shows that the hypoellipticity property of an operator does not imply the local solvability of the operator itself. However, given a hypoelliptic partial differential operator its adjoint is always locally solvable. In fact the adjoint of the Kannai operator, even though it has the same changing sign property of the principal symbol around the same set $S=\{x\in\mathbb{R}^n;x_1=0\}$, is locally solvable at each point of $\mathbb{R}^n$ (thus, also on $S$) due to the hypoellipticity property of the Kannai operator. Note also that the Kannai operator and its adjoint have the same principal symbol but different solvability property; therefore, unlike the principal type case, in presence of multiple characteristics the analysis of the principal symbol only is not sufficient to get conditions giving the local solvability and lower order terms are crucial.

Of course, we have locally solvable partial differential operators which are not the adjoints of hypoelliptic ones, therefore we can not get the local solvability from the hypoellipticity property. Moreover the problem of determining when a general partial differential operator is hypoelliptic is widely still open. 

These observations explain the reason why it is interesting to find classes of operators having the changing sign property of the principal symbol which are still locally solvable at the set around which this property is fulfilled.

A first class of operators having the changing sign property of the principal symbol generalizing the Kannai operator and its adjoint was studied by Colombini, Cordaro and Pernazza in \cite{CCP}, where they proved some $L^2$-local solvability results at the points where the property occurs by requiring a kind of condition $(\psi)$.

The class \eqref{Pintro}, introduced with A. Parmeggiani in \cite{FP}, represents a generalization of the class in \cite{CCP}. 
The local solvability property of operators of the form \eqref{Pintro} is studied both on $S$ and off the set $S$. In addition the problem on $S$ is also analyzed in the \textit{complex} case, where, in our notation, \textit{complex} case means that all the vector fields $iX_j$, $1\leq j\leq N$, appearing in the second order part of the operator are complex, while $iX_0$ is still \textit{real\textit{}}. Of course the \textit{complex} case is a generalization of the \textit{real} one; however the proof in the two cases is slightly different, and the technique used in the real case is relevant for a future generalization of the results in a pseudo-differential setting (see \cite{FP}). 

Recall once more that the class \eqref{Pintro} contains multiple characteristics partial differential operators for which the local solvability problem is very difficult to analyze  because of the complexity of the underlying geometry. Here the important results obtained for principal type pseudo-differential operators, as, in particular, the resolution of the Nierenberg-Treves conjecture due to Dencker \cite{D}, are not suitable, since they are based on the study of the principal symbol only, whereas, in presence of multiple characteristics, the lower order part of the symbol can not be neglected.

In the multiple characteristics setting we have solvability results by Mendoza \cite{M}, Mendoza and Uhlmann \cite{MU}, M\"uller and Ricci \cite{MR1}, Peloso and Ricci \cite{PeR}, Tr\'eves \cite{T}, and J.J. Kohn \cite{K} (complex vector fields). We also have recent work by Parenti and Parmeggiani \cite{PP} concerning the semi-global solvability in presence of multiple transverse symplectic characteristics.
There are also results by Beals and Fefferman \cite{BF} (see also Zuily \cite{Z} and Akamatsu \cite{A}) in which they give some conditions for the hypoellipticity of operators of the form $P^*$ where $P$ is contained in the class \eqref{Pintro}. Of course, their results for $P^*$ agrees with our local solvability results for $P$ in the class \eqref{Pintro}.

Beyond the class \eqref{Pintro} we shall study the local solvability problem for two classes of second order linear degenerate partial differential operators with non-smooth coefficients. These classes, analyzed by the author in \cite{F}, represent a variation with non-smooth coefficients of the class \eqref{Pintro}. Here the local solvability problem is considered only at points of $S$ (unlike the smooth case), since, far from $S$, these classes are contained in the main class \eqref{Pintro} and the results true for \eqref{Pintro} apply. 

All the results obtained for the classes presented above are proved by using the technique of a priori estimates, that is, the local solvability property follows directly from the validity of some estimates satisfied by the operators. In particular, depending on the kind of estimate satisfied by the operators, we gain a different kind of local solvability, that is, $H^s$ to $H^{s'}$ local solvability (see Definition \ref{local solvability s-s'}) with $s=-1,-1/2, 0,$ and $s'=0$, where the value of $s$ depends on the a priori inequality we are able to prove.

\section{Local solvability results for the main class}
We start with the analysis of the main class given by
\begin{equation}
\label{eqOperator}
P=\sum_{j=1}^NX_j^*f^pX_j+iX_0+a_0,
\end{equation}
where $p\geq 1$ is an integer, $X_j(x,D)$, $0\leq j\leq N$ ($D=-i\partial$), are \textit{homogeneous first order differential operators} 
(i.e. with no lower order terms) with smooth coefficients on an open set $\Omega\subset\mathbb{R}^n$
and with a \textit{real} principal symbol (in other words, the $iX_j(x,D)$ are \textit{real} vector fields; however, we shall also consider the case of $iX_1,\ldots,iX_N$ being \textit{complex}), 
$a_0$ is a smooth possibily complex-valued function, and $f\colon\Omega\longrightarrow\mathbb{R}$ a smooth function with $f^{-1}(0)\not=\emptyset$ and $df\bigl|_{f^{-1}(0)}\not=0$. The other classes we are going to treat later will be a variation of that introduced here.

As mentioned in the Introduction this class was inspired by a work of Colombini, Cordaro and Pernazza \cite{CCP} in which they study the local solvability of a class generalizing the Kannai operator and whose expression is given by
\begin{equation}
\label{CCP}
P=X^*(x,D)f(x)X(x,D)+iX_0(x,D)+a_0,
\end{equation}
where $iX$, $iX_0$ are real vector fields defined over an open set $\Omega\subset\mathbb{R}^n$, $X_0(x,D)\neq 0$ on $\Omega$, $f:\, \Omega\longrightarrow \mathbb{R}$ is \textit{analytic} and $a_0\in C^\infty(\Omega)$. They studied the local solvability of $P$ at points of $S=f^{-1}(0)$ around which the principal symbol changes sign, and proved, by requiring a kind of \textit{transversality} condition of $iX_0$ to the set $S$ that they called condition $(\psi)$, that the operator is locally solvable at each point of $S$.

The class \eqref{eqOperator} is a generalization of \eqref{CCP} since we suppose that $f$ is $C^\infty(\Omega)$ and not analytic, we consider a sum of terms in the second order part of the operator instead of a single term, and, in addition, we also consider the local solvability problem on $S$ in the \textit{complex} case, that is, when all the vector fields $iX_j$, $1\leq j\leq N$, are complex.

Let us remark that in \cite{FP} the local solvability of \eqref{eqOperator} is analyzed in several cases:
\begin{itemize}
  \item Solvability near $S=f^{-1}(0)$ in the \textit{real} case (that is when all the vector fields $iX_j$ are real) when:
  \begin{itemize}
  \item $p=2k+1$, $k\in\mathbb{Z}_+$;
   \item $p=2k$, $k\in\mathbb{Z}_+$;
\end{itemize}
  \item Solvability off $S=f^{-1}(0)$ when $p\in \mathbb{Z}_+$, in the \textit{real} case.
  \item Solvability near $S=f^{-1}(0)$ when $p=2k+1$, $k\in\mathbb{Z}_+$ in the  general \textit{complex} case, that is when all the vector fields $iX_j$, with $j\neq 0$, are complex (generalization of the \textit{real} case)
\end{itemize}

Let us also remark that when the exponent $p$ is even we do not have the changing sign property of the principal symbol; however, since the operator is still highly degenerate in this case, the local solvability is not guaranteed.
On the other hand, as we do not ask for the $iX_j$, $1\leq j\leq N$, non-degeneracy conditions, the local solvability far from $S$ has to be studied as well.

Depending on the hypotheses satisfied by the class at a point $x_0\in \Omega$ (or, more precisely, in a neighborhood of $x_0$) we can get a different kind of local solvability at $x_0$. We recall below the definition of local solvability that we shall use throughout the paper.

\begin{definition}[Local solvability in the sense $H^s$ to $H^{s'}$]\label{local solvability s-s'}
We say that P is $H^s$ to $H^{s'}$ locally solvable at $x_0\in\Omega$ if there exists a compact $K\subset \Omega$, with $x_0\in \mathring{K}=U$ (where $ \mathring{K}$ is the interior of $K$), such that for all $f\in H^s_{\mathrm{loc}}(\Omega)$ there is $u\in H^{s'}_{\mathrm{loc}}(\Omega)$ which solves $Pu=f$ in $U$. 
\end{definition}
In order to get solvability results for \eqref{eqOperator} we need to impose some conditions, which, in particular, are given in terms of hypotheses on the vector fields involved in the expression of the operator. Since we are going to analyze the \textit{real} case first, we start by giving the hypotheses in this case. However, in the \textit{complex} case, we will assume a suitable generalization of the hypotheses required in the real setting.
We shall denote by $H_j$ the Hamilton vector field associated with the symbol  $X_j(x,\xi)$ of $X_j(x,D)$, by $S=f^{-1}(0)$ the zero set of the function $f$ involved in the expression of the operator and by $\Sigma=\bigcap_{j=0}^N\Sigma_j$, where $\Sigma_j$ is the characteristic set of the operator $X_j(x,D)$ ($X_j$ for short), namely $\Sigma_j=\{(x,\xi)\in\Omega\times\mathbb{R}^n| X_j(x,\xi)=0\}$.

We shall call (H1), (H2) and (H3) the following conditions:
\begin{itemize}
  \item[(H1)] $iX_0(x,D)f>0$ on $S$;
  \item[(H2)] $\forall j=1,...,N$, $\forall K\subset\Omega$, $\exists C>0$ such that
  \begin{equation*}
\{X_j,X_0\}(x,\xi)^2\leq C\sum_{k=0}^N X_k(x,\xi)^2,\quad \forall(x,\xi)\in K\times \mathbb{R}^n;
\end{equation*}
  \item[(H3)] let $\rho\in\Sigma,\,V(\rho)=\mathrm{Span}\{H_0(\rho),...,H_N(\rho)\}$, $J=J(\rho)$$\subset\{0,...,N\}$ be such that $\{H_j(\rho)\}_{j\in J}$ is a basis of $V(\rho)$, and $M(\rho)=$$\left[\{X_j,X_{j'}\}(\rho)\right]_{j,j'\in J}$. We say that (H3) is satisfied at $x_0$ if $\pi^{-1}(x_0)\cap\Sigma\neq\emptyset$ and
  \begin{equation*}
\mathrm{rank}\,M(\rho)\geq 2,\quad \forall \rho\in \pi^{-1}(x_0)\cap\Sigma.
\end{equation*} \end{itemize} 
Here $\{\cdot,\cdot\}$ denotes the Poisson bracket with respect to the standard symplectic from $\sigma$ on $T^*\Omega$, $\pi^{-1}\{x_0\}$ is the fiber of $x_0$ and $\{X_j,X_0\}(x,\xi)$, $X_j(x,\xi)$ are the total (principal because of homegeneity) symbols of $i[X_j,X_0]$ and $X_j$ respectively. 

\begin{remark}
It is possible to prove that condition (H3) at a point $x_0$ is equivalent to the H\"ormander condition at step 2 for the system of vector fields $\{iX_j\}_{j=0,...,N}$ at the point $x_0$,  and thus also in a sufficiently small neighborhood of $x_0$.
Note also that when $\pi^{-1}(x_0)\cap\Sigma=\emptyset$ then the system of vector fields $\{X_j\}_{j=0,...,N}$ is elliptic in a neighborhood of $x_0$.
\end{remark} 

We now give the statement of the first result concerning the local solvability of $P$ of the form \eqref{eqOperator} around the set $S=f^{-1}(0)$ when $p$ is an odd  integer in the real case. 

\begin{theorem}[Solvability near $S$ when $p=2k+1$ in the real case]
\label{thmLocSolv}
Let the operator $P$ in \eqref{eqOperator} satisfy hypotheses (H1) and (H2). 
\begin{itemize}
\item[(i)] Let $k=0$. Then for all $x_0\in S$ with $\pi^{-1}(x_0)\cap\Sigma\not=\emptyset$ and at which 
hypothesis (H3) is fulfilled,
there exists a compact $K\subset W$ with $x_0\in U=\mathring{K}$ such that
for all $v\in H^{-1/2}_{\mathrm{loc}}(\Omega)$ there exists $u\in L^2(\Omega)$ solving $Pu=v$ in $U$
(hence we have $H^{-1/2}$ to $L^2$ local solvability). 
\item[(ii)] Let $k=0$.
Then, for all $x_0\in S$ for which $\pi^{-1}(x_0)\cap\Sigma=\emptyset$
there exists a compact $K\subset\Omega$ with $x_0\in U=\mathring{K}$ such that
for all $v\in H^{-1}_{\mathrm{loc}}(\Omega)$ there exists $u\in L^2(\Omega)$ solving $Pu=v$ in $U$
(hence we have $H^{-1}$ to $L^2$ local solvability). 
\item[(iii)] If $k\geq 1$ and $x_0$ is any given point of $S$,
or $k=0$ and $x_0\in S$ (with $\pi^{-1}(x_0)\cap\Sigma\not=\emptyset$) is such that (H3) is not satisfied at $x_0$,
then there exists a compact  $K\subset\Omega$ with $x_0\in U=\mathring{K}$ such that
for all $v\in L^2_{\mathrm{loc}}(\Omega)$ there exists $u\in L^2(\Omega)$ solving $Pu=v$ in $U$
(hence we have $L^2$ to $L^2$ local solvability). 
\end{itemize}

\end{theorem}

\begin{remark}
The expression of the class \eqref{eqOperator} under consideration is invariant under change of variables as well as conditions (H1), (H2) and (H3). These properties make our result general in the sense that, given an operator in the class satisfying the hypotheses of the theorem, after a change of variables the new operator (that will have the same form) satisfies the same hypotheses and the theorem still applies.
\end{remark}

As one can see from the statement of Theorem \ref{thmLocSolv}, depending on the geometric relations among the vector fields $iX_j$, $1\leq j\leq N$, we can have different kinds of local solvability for $P$ of the form \eqref{eqOperator}. In particular, when $k=0$, we can have points of $S$ where we have a kind of local solvability and others where we have a different kind of local solvability. This property is shown by the following example.

\begin{example}\label{ex1}
Let $\Omega\subset\mathbb{R}^3$ be such that $\{x\in\mathbb{R}^3; x_1=-1\}\subset \Omega$, and write $\Omega$ as $\Omega=\Omega_-\cup \{x\in\mathbb{R}^3;x_1=-1\}\cup\Omega_+$, with $\Omega_{\pm}=\{x\in\Omega;x_1\gtrless -1\}$.
We now take $f(x)=x_2+x_2^3/3-x_1x_3$, $X_0(x,D)=D_2-x_1D_3$, $X_1(x,D)=D_1-x_3D_3$ and $X_2(x,D)=(1+x_1)D_3$,  and study the local solvability in $S=f^{-1}(0)$ of the operator (defined over $\Omega$)
$$P=\sum_{j=1}^2X_j^*fX_j+iX_0+a_0\quad (k=0).$$
$P$ is an operator in \eqref{eqOperator} and conditions (H1) and (H2) are satisfied, since $iX_0f=1+x_2^2+x_1^2$ and 
$$\{X_1,X_0\}=-X_2,\quad \{X_2,X_0\}=(2+x_1)\xi_3,\quad \{X_1,X_2\}=0.$$
Moreover $\Sigma=\{(x,\xi)\in\Omega\times(\mathbb{R}^3\setminus\{0\});x_1=-1,\xi_1=x_3\xi_3,\xi_2=-\xi_3,\xi_3\neq 0\}$, therefore $\Sigma\cap\pi^{-1}(\Omega_\pm)=\emptyset$ and $\Sigma\cap \pi^{-1}(\{x_1=-1\})\neq\emptyset$.
Now we analyze the local solvability of $P$ in $S$ in the cases $x_0\in S\cap \Omega_\pm$ and $x_0\in S\cap\{x_1=-1\}$.\\
For all $x_0\in S\cap \Omega_\pm$ we have $\pi^{-1}(x_0)\cap\Sigma=\emptyset$. Thus, by point (ii) of the theorem, for all $x_0\in S\cap \Omega_\pm$, $P$ is $H^{-1}$ to $L^2$ locally solvable at $x_0$.\\
For all $x_0\in S\cap\{x_1=-1\}$ we have $\pi^{-1}(x_0)\cap\Sigma\neq\emptyset$. Moreover $H_{0}$, $H_{1}$ and $H_{2}$ are linearly independent, and condition (H3) is satisfied at $x_0$ for all $x_0\in S\cap\{x_1=-1\}$. It follows, by point (i) of Theorem \ref{thmLocSolv}, that for all $x_0\in S\cap\{x_1=-1\}$, $P$ is $H^{-1/2}$ to $L^2$ locally solvable at $x_0$.
\end{example}

Another relevant observation is that the class \eqref{eqOperator} contains operators being not the adjoints of hypoelliptic ones. This is important to remark, since, otherwise, the local solvability property of an operator in the class would follow directly by the hypoellipticity property of its adjoint. Example \ref{ex2} below shows that the class \eqref{eqOperator} contains in fact operators whose adjoint is not hypoelliptic.

\begin{example}\label{ex2}
Let $g\in C_0^\infty(\mathbb{R};\mathbb{R})$ be such that $g\not\equiv 0$ and $g(0)\neq 0$. Following a result by Zuily \cite{Z} we know that the operator
$$L_k(x,D)=iD_2-(x_2-g(x_1))^kD_1^2,\quad (x_1,x_2)=x\in\mathbb{R}^2,$$
is hypoelliptic for all $k\neq 1$ (thus $L_{k\neq 1} ^*$ is locally solvable).\\
We now consider the operator
$$P=X_1^*fX_1 +iX_0+a_0,$$
where $X_1=D_1$, $X_0=D_2+g'(x_1)D_1$, $f(x)=x_2-g(x_1)$, and $a_0(x)=g''(x_1)$. Note that $P$ is of the form \eqref{eqOperator} and that $P^*$ is not hypoelliptic since $P^*=-L_1$. Moreover, by requiring $|g'(x_1)|\leq c<1$ (so that $P$ satisfies condition (H1)), by the previous theorem we have that $P$ is $H^{-1}$ to $L^2$ locally solvable in $S$.
\end{example}

\subsection*{Sketch of the proof of Theorem \ref{thmLocSolv}}
By functional analysis arguments it is possible to prove that the local solvability property of a partial differential operator is equivalent to the validity of an a priori estimate satisfied by the operator.
We then rephrase the definition of $H^s$ to $H^{s'}$ local solvability in terms of a priori estimates in the following way.
\begin{definition}[Definition of local solvability via a priori estimates]\label{Def Solvability}
A partial differential operator $P$, defined over an open set $\Omega\subset\mathbb{R}^n$, is $H^s$ to $H^{s'}$ locally solvable at $x_0\in\Omega$ if there exists a compact $K\subset \Omega$, with $x_0\in \mathring{K}=U$, and a positive constant $C$ such that, for all $u\in C_0^\infty(U)$,
\begin{equation}\label{SolvabilityEst}
 C\n P^*u\n_{-s'}\geq \n u\n_{-s}.
\end{equation}
\end{definition}
We will often refer to \eqref{SolvabilityEst} as the \textit{solvability estimate}.

The first step of the proof consists in the following estimate satisfied by $P^*$ that we shall call \textit{main estimate}.

\begin{proposition}\textit{[Main Estimate}]\label{propP*}
Let $P$ be as in (\ref{eqOperator}) and satisfy (H1).
For all $x_0\in S$ there exists a compact $K_0\subset\Omega$, with $x_0\in\mathring{K}_0$, 
constants $c=c(K_0),C=C(K_0)>0$ and $\varepsilon_0=\varepsilon_0(K_0)$ with $\varepsilon_0(R)\to 0$ as
the compact $R\searrow\{x_0\}$, such that for all compact $K\subset K_0$ 
\begin{equation}
\hspace{-.8cm}
|\!|P^*u|\!|_0^2\geq \frac18|\!|X_0u|\!|_0^2+c\Bigl(\widehat{P}_0(x,D)u,u\Bigr)-C|\!|u|\!|_0^2,
\label{MainEst}\end{equation}
for all $u\in C_0^\infty(K)$, where
\begin{equation}\label{P0}
\widehat{P}_0=X_0^*X_0+\sum_{j=1}^N(X_j^*f^{2k}X_j-\varepsilon_0^2[X_j,X_0]^*f^{2k}[X_j,X_0]),
\end{equation}
and $(\cdot,\cdot)$ is the $L^2$-scalar product.
\end{proposition}
For the proof of the previous result see \cite{FP}.

Following Definition \ref{Def Solvability} we have that, in order to prove point (i), (ii) and (iii) of Theorem \ref{thmLocSolv} we have to show that for all $x_0\in S$ there exists a neighborhood $U$ of $x_0$ such that the operator $P^*$, which is the adjoint of $P$ of the form \eqref{eqOperator}, satisfies the estimate \eqref{SolvabilityEst} with $(s,s')=(-1/2,0)$, $(s,s')=(-1,0)$ and $(s,s')=(0,0)$ respectively. 
This is possible starting from the main estimate by estimating the term $\Bigl(\widehat{P}_0(x,D)u,u\Bigr)$ from below and then absorbing the $L^2$ error by using a Poincar\'e inequality on $X_0$ (which is nondegenerate near $S$ by condition (H1)). In particular we get
\begin{itemize}
  \item $H^{-1/2}$ to $L^2$ local solvability at $x_0\in S$ (point (i) of the theorem), by using the \textit{Melin inequality} on $\widehat{P}_0(x,D)$;
  \item $H^{-1}$ to $L^2$ local solvability at $x_0\in S$ (point (ii) of the theorem), by using the \textit{G\aa rding inequality} on $\widehat{P}_0(x,D)$;
  \item $L^2$ to $L^2$ local solvability at $x_0\in S$ (point (iii) of the theorem), by using the \textit{Fefferman-Phong inequality} on $\widehat{P}_0(x,D)$.
\end{itemize}

Note that in the expression of $\widehat{P}_0$ the constant $\varepsilon_0^2$ is a positive constant depending on the compact set $K_0$ that shrinks as $K_0$ is shrunk around $x_0$. The latter property and hypotheses (H2) and (H3) are very important in order to get a control of the commutators in \eqref{P0} and to apply the Melin, respectively the G\aa rding and the Fefferman-Phong inequality on $\widehat{P}_0$ (see \cite{FP}).

\subsection*{Other results}
As mentioned before the local solvability of the class \eqref{eqOperator} is also proved on $S$ when $p=2k$, $k\geq 1$, and off $S$ when $p\geq 1$, $p\in \mathbb{Z}$. In these cases we do not have the changing sign property of the principal symbol anymore. However the operator is still degenerate and therefore the local solvability is not guaranteed. Moreover the results above can be proved under weaker assumptions. In fact hypotheses (H1) and (H2), that was fundamental in the previous case, are not assumed  now, and the proof in both the cases described here is based on the use of Carleman estimates. We give below the statement of these results.

\begin{theorem}[Solvability near $S$, $p=2k$, $k\geq 1$]
Let $P$ of the form \eqref{eqOperator} be such that $iX_0f\neq 0$ on $S$. Then for all $x_0\in S$ there exists a compact $K\subset \Omega$, with $x_0\in\mathring{K}=U$, such that $P$ is $L^2$ to $L^2$ locally solvable in $U$.
\end{theorem}
\begin{theorem}[Solvability off $S$, $p\in \mathbb{Z}_+$]\label{offS}
Let $P$ of the form \eqref{eqOperator} be such that for every $x_0\in \Omega\setminus S$ there exists $j_0$, $1\leq j_0\leq N$, for which $X_{j_0}$ is nonsingular at $x_0$. Then $P$ is $L^2$ to $L^2$ locally solvable in $\Omega\setminus S$.
\end{theorem}
For the proof of these theorems see \cite{FP}.
\begin{remark}
In the cases covered by the preceding theorems we can only prove $L^2$ to $L^2$ local solvability results. However, we do not impose specific conditions on the system of vector fields $\{iX_j\}_{j=0,...,N}$, therefore the setting we can consider here is more general, since the changing sign property of the principal symbol does not occur. Of course, by adding some additional conditions on the system $\{X_j\}_{j=1,...,N}$, at least off the set $S$, one can prove better solvability results for the class \eqref{eqOperator}.

Note that in the odd exponent case we have better local solvability results at a point $x_0\in S$ only when the exponent $p=1$, that is, the degeneracy carried by the function $f$ is 1, and the vector fields $\{X_j\}_{j=0,...,N}$ satisfy additional conditions (that is, when the system satisfies the H\"{o}rmander condition at step 2 at $x_0$ or is elliptic at $x_0$). In general, when the degeneracy carried by the function $f^p$ is greater than or equal to 2, that is, $p\geq 2$, our results give $L^2$ to $L^2$ local solvability in $S$. This observation is to remark that the degree of  degeneracy of the operator deeply influences its regularity properties.

As regards the problem far from $S$, where the degeneracy carried by the function $f^p$ is zero, that is, the degeneracy of $P$ can only follow by the degeneracy of the system of vector fields $\{iX_j\}_{j=1,...,N}$, then, again, by asking the system $\{iX_j\}_{j=0,...,N}$ to fulfill additional conditions as before, one can get better results. 

\end{remark}
\subsection*{The complex case}
It is possible to generalize Theorem \ref{thmLocSolv} to a \textit{complex} case, that is, when all the vector fields $iX_j$, $1\leq j\leq N$, in \eqref{eqOperator} are complex. Here the hypotheses are a suitable generalization of those given in the real case. Besides (H1) we shall assume:
\begin{itemize}
\item[(H2$'$)] \textit{For all $j=1,\ldots,N$ and for all compact $K\subset\Omega$ there exists $C_{K,j}>0$ such that}
$$|\{X_j,X_0\}(x,\xi)|^2\leq C_{K,j}\sum_{j'=0}^N|X_{j'}(x,\xi)|^2,\,\,\,\,\forall(x,\xi)\in \pi^{-1}(K);$$
\item[(H3$'$)] \textit{For all compact $K\subset\Omega$ there exists $C_K>0$ such that}
$$|\sum_{j=1}^N\{\bar{X}_j,X_j\}(x,\xi)|^2\leq C_K\sum_{j=0}^N|X_j(x,\xi)|^2,\,\,\,\,\forall(x,\xi)\in \pi^{-1}(K);$$
\item[(H4$'$)] \textit{For all compact $K\subset\Omega$ there exists $C_K>0$ such that}
$$|\sum_{j=1}^N\{\overline{\{X_j,X_0\}},\{X_j,X_0\}\}(x,\xi)|^2\leq C_K\sum_{j=0}^N|X_j(x,\xi)|^2,\,\,\,\,\forall(x,\xi)\in \pi^{-1}(K).$$
\end{itemize}

These new suitable conditions, as in the real case, permit the control of the commutator part inside the operator $\widehat{P}_0$ appearing in the main estimate \eqref{MainEst}, which holds in this complex case as well, and the application of the H\"{o}rmander and the Fefferman-Phong inequality for  $\widehat{P}_0$.

The theorem true in this setting is the following.

\begin{theorem}\label{thmComplex}
Let $P$ be given as in (\ref{eqOperator}) satisfying hypotheses (H1), (H$2'$), (H$3'$) and (H$4'$) above.
\begin{itemize}
\item [(i)] If $k\geq 0,$ then for all $x_0\in S$ there exists a compact $K\subset\Omega$ with $x_0\in\mathring{K}=U$
such that $L^2$ to $L^2$ local solvability holds on $U$.
\item [(ii)] If $k=0,$ suppose in addition that
there exists $r\geq 1$ such that for all $x_0\in S$ for which $\pi^{-1}(x_0)\cap\Sigma\not=\emptyset$ 
there exists a neighborhood $W_{x_0}\subset\Omega$ of $x_0$ such that
\begin{equation}
\dim\mathscr{L}_r(x)=n,\,\,\,\,\forall x\in W_{x_0}.
\label{eqHor}\end{equation}
Then for each $x_0\in S$ there exists a compact $K\subset\Omega$ with $x_0\in\mathring{K}=U$ such that we
have $H^{-1/r}$ to $L^2$ solvability on $U$.
\end{itemize}
\end{theorem}
Here $\mathscr{L}_r(x)$ stands for the Lie algebra generated by the real and the imaginary parts of the vector fields $iX_j$, $0\leq j\leq N$, and by their commutator up to length $r$.

The proof of this theorem follows the same line of the proof of Theorem \ref{thmLocSolv}. One starts from the main estimate \eqref{MainEst} to get the solvability estimate by controlling the term $(\widehat{P}_0u,u)$ and cancelling the $L^2$ error. Here, since the vector fields are not real, we estimate $\widehat{P}_0$ by using different methods, that is, in particular, depending on the hypotheses, we can use a combination of the subelliptic H\"ormander inequality on the system of vector fields $\{iX_j\}_{j=0,...,N}$ satisfying the H\"ormander condition at step $r\geq 1$ and the Fefferman-Phong inequality. Finally, once more, we absorb the $L^2$ error by using a Poincar\'e inequality on $X_0$ and obtain the suitable solvability estimate. For more details see \cite{FP}.

To conclude this section, we give an example of an operator in the class.
\begin{example}
Let $\omega\in C^\infty(\mathbb{R}^2_{x_1,x_2};\mathbb{R})$  be such that $\partial_{x_2}\omega\neq 0$, and consider $f(x)=x_4+g(x_1,x_2,x_3)$ ($x\in\mathbb{R}^4$), $g\in C^\infty(\mathbb{R}^3_{x_1,x_2,x_3};\mathbb{R})$, and $X_0,X_1,X_2$ equal to
$$Z_0(x,D)=D_4,\quad Z_1(x,D)=D_1+ix_2^kD_4,\quad Z_2(x,D)=e^{i\omega(x_1,x_2)}D_2.$$
Then, by the previous theorem, we have that
$$P=\sum_{j=1}^2X_j^*fX_j+iX_0$$
is $H^{-\frac{1}{k+1}}$ to $L^2$ locally solvable at $S$.
\end{example}
\section{The non-smooth coefficients case}
We next discuss the local solvability of two classes of second order linear degenerate partial differential operators with non-smooth coefficients which are a variation of the main class \eqref{eqOperator}.
\subsection*{The first class with non-smooth coefficients}
The first class we shall analyze is given by
\begin{equation}
\label{GP}
P=\sum_{j=1}^NX_j^*f|f|X_j+iX_0+a_0,
\end{equation}
where $X_j=X_j(D)$, $1\leq j\leq N$, are homogeneous first order differential operators (in other words $iX_j$, $0\leq j\leq N$ are vector fields) with real or complex constant coefficients (the two cases will be analyzed separately), $iX_0=iX_0(x,D)$ is a real vector field with \textit{affine} coefficients, $f\not\equiv \text{constant}$ is an affine function, and $a_0$ is a continuous  function on $\mathbb{R}^n$ with complex values.
Note that $P$ of the form \eqref{GP} could have $C^{0,1}$ or $C^{1,1}$ coefficients depending on the transversality or tangency of the vector fields $iX_j$, $j\neq 0$, to the set $f^{-1}(0)=S$. Moreover the class has an important invariance property, that is, it is invariant under \textit{affine} change of variables.

Since the operator is degenerate on $S$,  and can be degenerate far from $S$ as well, as before the local solvability is not guaranteed. 
The local solvability of $P$ is studied in the neighborhood of $S$ where the principal symbol changes sign, namely, the principal symbol changes sign in the neighborhood of each point of $\pi^{-1}(S)$. As regards the solvability far from $S$, we do not analyze this problem here, since, in this case, the class \eqref{GP} is contained in the class \eqref{eqOperator} and Theorem \ref{offS} still applies.
We shall consider the problem for the class \eqref{GP} in two cases, that is, in the \textit{real} case when all the vector fields $iX_j$, $1\leq j\leq N$, in the second order part of the operator are real, and also in the \textit{complex} case when all the vector fields $iX_j$, $1\leq j\leq N$ are complex. Note that the vector field $iX_0$ is always taken with real coefficients. Let us remark that in the complex case one needs an additional condition in order to prove the result. 

Since one deals with operators with non-smooth coefficients, one is able to prove an $L^2$-local solvability result at the points of $S$ with no possibility to gain regularity. In particular, due to the non-smoothness of $f|f|$, one needs to clarify the meaning of $L^2$-local solvability (for more information about solvability see \cite{Ho4} and \cite{L}).

\begin{definition}\label{DefSolvability2}
Given a partial differential operator $P$, defined on an open set $\Omega\subseteq \mathbb{R}^n$, such that both $P$ and its adjoint $P^*$ have at least $L^\infty$ coefficients, one says that $P$ is $L^2$-locally solvable in $\Omega$ if for any given $x_0\in\Omega$ there is a compact set $K\subset \Omega$ with $x_0\in U=\mathring{K}$, such that for all $f\in L^2_{\mathrm{loc}}(\Omega)$ there exists $u\in L^2(U)$ such that 
$$(u, P^*\varphi)=(f,\varphi),\quad  \forall \varphi\in C_0^\infty(U),$$
 where $(\cdot,\cdot)$ is the $L^2$ inner product.
 \end{definition}

As regards the hypotheses on the class, we assume the following formulation in this setting, analogous to the previous ones:
\begin{itemize}
\item[(H1)] $iX_0f(x)>0$ for all $x\in S:=f^{-1}(0)$;
\item[(H2)] for all  $1\leq j\leq N$ there exists a constant $C>0$ such that \\
$$|\{X_j,X_0\}(\xi)|^2\leq C\sum_{j=1}^N|X_j(\xi)|^2,\,\,\,\,\forall \xi\in \mathbb{R}^n;$$
\item[(H3)] when $X_1,..., X_N$ have complex coefficients, then $X_jf=0$, $\forall j=1,...,N.$
\end{itemize}

Condition (H3) is used only in the complex case and means that each vector field $X_j$, with $j\neq 0$, is tangent to $S=f^{-1}(0)$.

We now are in a position to give two results concerning the solvability in the real and in the complex case respectively.

\begin{theorem}[Real case]\label{thmLocSolv2}
Let $P$ be of the form (\ref{GP}) such that all the vector fields $iX_j$ have real coefficients and hypotheses (H1),(H2) are satisfied, and let $S$ be the zero set of $f$. Then for all $x_0\in S$ there exists a compact $K\subset \mathbb{R}^n$ with $U=\mathring{K}$ and $x_0\in U$, such that for all $v\in L^2_{\mathrm{loc}}(\mathbb{R}^n)$ there exists $u\in L^2(U)$ solving $Pu=v$ in $U$ (in the sense of Definition \ref{DefSolvability2}).
\end{theorem}

\begin{theorem}[Complex case]\label{solvcomplex}
Let $P$ be of the form (\ref{GP}) such that hypotheses (H1) to (H3) are satisfied, and let $S$ be the zero set of $f$. Then for all $x_0\in S$ there exist a compact $K\subset \mathbb{R}^n$ with $U=\mathring{K}$ and $x_0\in U$, such that for all $v\in L^2_{\mathrm{loc}}(\mathbb{R}^n)$ there exists $u\in L^2(U)$ solving $Pu=v$ in $U$ (in the sense of Definition \ref{DefSolvability2}).
\end{theorem}

Both the results above are proved by using the technique of a priori estimates following the scheme of the the proof of Theorem \ref{thmLocSolv}. One starts with the proof of an estimate satisfied by $P^*$ involving a suitable operator $\widehat{P}_0$ having non-smooth coefficients. Thanks to the hypotheses one can use Fourier analysis to control $\widehat{P}_0$, and finally, once more by using a Poincar\'e inequality on $X_0$, cancel the $L^2$ error to get the solvability estimate giving the result.

For the proof of these results see \cite{F}.

\begin{remark}
Note that in the complex case we need to require a tangency condition on the vector fields $iX_j$, $1\leq j\leq N$. This condition can be removed in a special case, that is, when $N=1$ and $X_1=\langle \alpha, D\rangle$, $\alpha\in\mathbb{C}^n$, is such that
$$\mathrm{dim}\big(\mathrm{Span}_{\mathbb{R}}\{\mathsf{Re}(\alpha),\mathsf{Im}(\alpha)\}\big)=1.$$
Observe that the  condition above is equivalent to say that $X_1$ is essentially real, namely $\overline{X_1(\xi)}=\zeta X_1(\xi)$, $\zeta\in\mathbb{S}^1\subset\mathbb{C}$.
With these assumptions, still under hypotheses (H1) and (H2), $\forall x_0\in S$, there exists a compact $K\subset\mathbb{R}^n$, with $x_0\in\mathring{K}=U$, on which $P$ is $L^2$ to $L^2$ locally solvable in $U$.
\end{remark}

\subsection*{The second class with non-smooth coefficients}
We present here the second class with non-smooth coefficients whose expression is given by
\begin{equation}\label{PP}
P=\sum_{j=1}^NX_j^*|f|X_j+iX_0+a_0,
\end{equation}
where $X_j=X_j(x,D)$, $0\leq j\leq N$, are homogeneous first order differential operators with smooth coefficients defined on an open set $\Omega\subseteq\mathbb{R}^n$, $f: \Omega\longrightarrow\mathbb{R}$ is a $C^1$ function with $S=f^{-1}(0)\neq \emptyset$ and $df|_S\neq 0$, and $a_0$ is a continuous possibly complex valued function.
Now the vector fields $iX_j$, $1\leq j\leq N$, have smooth real or complex coefficients  (we \textbf{do not} assume constant coefficients anymore), $f$ is not affine but a smooth real function, and $iX_0$ has smooth real variable coefficients.

Once more we study the local solvability of \eqref{PP} at points of $S$, since, off the set $S$, again, this class is contained in the main class \eqref{eqOperator}.
Let us also remark that we do not have the changing sign property of the principal symbol around $S$. However, the operator is still degenerate and characterized by a less regularity of the coefficients.
In fact an operator of the form \eqref{PP} could have $C^{0,1}$ or $L^\infty$ coefficients depending on the tangency or transversality of the vector fields $iX_j$, $1\leq j\leq N$, to the zero set of $f$.

We now assume just hypothesis (H1), that is $iX_0(x,D)f>0$ on $S$, and get the following $L^2$-solvability result.
\begin{theorem}
Let the operator $P$ in (\ref{PP}) satisfy hypothesis $\mathrm{(H1)}$. Then for all  $x_0\in S$ there exists a compact $K\subset\Omega$ with $x_0\in U=\mathring{K}$, such that for all $v\in L^2_{\mathrm{loc}}(\mathbb{R}^n)$ there exists $u\in L^2(U)$ solving $Pu=v$ in $U$ (in the sense of Definition \ref{DefSolvability2}).
\end{theorem}

The proof of this theorem consists in the proof of the solvability estimate with $(s,s')=(0,0)$ starting from a Carleman estimate on $P^*$. For the entire  proof see \cite{F}.

\section{Final remarks}
For the main class, which is the one with smooth coefficients, we proved a result in a complex case, that is, when all the vector fields $iX_j$, $1\leq j\leq N$, in \eqref{eqOperator} are complex. However, we asked for the vector field $iX_0$ to be real, therefore the next goal is to find a general result for the class \eqref{eqOperator} when all the vector fields involved in the expression of the operator are complex. This is work in progress \cite{FP1}.

Another important aspect related to the main class is given by the regularity problem. As one can see by the results true for this class, we are able to state the existence of $L^2$ solutions for operators of the form \eqref{eqOperator}. Of course we expect, at least for some subclasses of \eqref{eqOperator}, to find solutions with higher regularity property, that is, given a source term $v$ in $H^s$, with $s>0$, there exists a solution $u$ in $H^{s+\varepsilon}$, with $\varepsilon>0$, of the problem $Pu=v$. In fact it is possible to prove that the adjoint of the Kannai operator, which, in particular, is contained in the main class \eqref{eqOperator}, is locally solvable with a loss of one derivative on $\mathbb{R}^n$ (thus also on $S$), namely, for all $x_0\in \mathbb{R}^n$ there is a neighborhood $U$ such that for all $v\in H^{-s+1}(\mathbb{R}^n)$ there exists a solution $u\in H^{-s}(U)$ of the problem $Pu=v$ in $U$ for all $s>-n/2$, where the condition on the index $s$ is related to the application of a Poincar\'e inequality in Sobolev spaces (see \cite{L}). To prove this property the hypoellipticity of the Kannai operator is strongly used, therefore a similar result for operators which are not adjoints of hypoelliptic ones is more difficult to obtain. We do not expect that all the operators in the class have this property; however, a classification of \eqref{eqOperator} into subclasses  depending on the regularity property is a very interesting argument to examine. At the moment the regularity problem for this class is open.

\bibliographystyle{alpha}

\end{document}